\documentclass[12pt]{article}

\usepackage[utf8]{inputenc}
\usepackage{xcolor}
\usepackage{hyperref}
\usepackage{amssymb,amsmath,amsfonts,graphicx,color,multicol,amsthm}
\usepackage{subcaption}
\usepackage[hmargin=.7in,vmargin=1.5in]{geometry}
\usepackage[T1]{fontenc}

\newtheorem{thm}{Theorem}[section]
\newtheorem{lem}[thm]{Lemma}

\newtheorem{rem}[thm]{Remark}

\begin{document}

\title{Linear multistep methods and global Richardson extrapolation}

\author{Imre Fekete\thanks{{\texttt{imre.fekete@ttk.elte.hu}}, corresponding author, Department of Applied Analysis and Computational Mathematics, ELTE E\"otv\"os Lor\'and University, and MTA-ELTE Numerical Analysis and Large Networks Research Group, Pázmány P. s. 1/c, H-1117 Budapest, Hungary}, Lajos L\'oczi\thanks{{\texttt{LLoczi@inf.elte.hu}}, Department of Numerical Analysis, ELTE E\"otv\"os Lor\'and University, Pázmány P. s. 1/c, H-1117 Budapest, Hungary, and Department of Differential Equations, BME Budapest University of Technology and Economics}}
\date{\today}

\maketitle

\begin{abstract}
In this work, we study the application the classical Richardson extrapolation (RE) technique to accelerate the convergence of sequences resulting from linear multistep methods (LMMs) for solving initial-value problems of systems of ordinary differential equations numerically. The advantage of the LMM-RE approach is that the combined method possesses  higher order and favorable linear stability properties in terms of $A$- or $A(\alpha)$-stability, and existing LMM codes can be used without any modification.
\end{abstract}


\section{Introduction}

Richardson extrapolation (RE) \cite{richardson1911, richardson1927} is a classical technique to accelerate the convergence of numerical sequences depending on a small parameter, by eliminating the lowest order error term(s) from the corresponding asymptotic expansion.
When the sequence is generated by a numerical method solving the initial-value problem 
\begin{equation}\label{ODE}
y'(t)=f(t,y(t)), \quad y(t_0)=y_0,
\end{equation}
 the parameter in RE can be chosen as the discretization step size $h>0$. The application of RE to sequences generated by one-step---e.g., Runge--Kutta---methods is described, for example, in \cite{butcher, hairernorsettwanner}. In \cite{zlatev}, global (also known as passive) or local (active) versions of RE are implemented with Runge--Kutta sequences. These combined methods can find applications in air pollution problems \cite{zlatev2022} or in machine learning \cite{falgout2021}, for example.


In this paper, we analyze the application of global Richardson extrapolation (GRE) to sequences generated by linear multistep methods (LMMs) approximating the solution of \eqref{ODE}.
We will refer to a $k$-step LMM as the underlying LMM, and its recursion has the usual form
\begin{equation}\label{multistepdefn}
\sum_{j=0}^k \alpha_j y_{n+j}=\sum_{j=0}^k h \beta_j f_{n+j},
\end{equation}
where $f_m:=f(t_m, y_m)$, and the numbers $\alpha_j\in\mathbb{R}$ and $\beta_j\in \mathbb{R}$ ($j=0, \ldots, k$) are the given method coefficients with $\alpha_k\ne 0$. The LMM is implicit, if $\beta_k\ne 0$.

In Section \ref{sec:main}, given an underlying LMM, we define its extrapolated version, referred to as LMM-GRE,  and investigate its convergence. Then, we carry out linear stability analysis for the LMM-GREs. In Section \ref{sec24}, we focus on the BDF family as underlying LMMs due to their good stability properties, although the results from Sections \ref{sec21}--\ref{sec23} are clearly applicable to other LMM families as well.
The numerical experiments in Section \ref{sec:numerics} demonstrate the expected convergence order---here, we use several different types of LMMs with GRE to solve \eqref{ODE}.


As a conclusion of this study, we see that\\
 (i) to implement a  LMM-GRE,  existing LMM codes can directly be used, thanks to the simple linear combination appearing in definition \eqref{GREdef}; moreover, \\
(ii) the higher computational cost of  a LMM-GRE is compensated by its higher convergence order and favorable linear stability properties.

\begin{rem} In \cite[Section 3.4]{zlatev}, the authors comment on the possible combination of LMMs and \emph{local} Richardson extrapolation. Working out the necessary details and convergence theorems for this case could be the subject of a future study.
\end{rem}

\section{Results for LMM-GREs}\label{sec:main}

\subsection{Definition}\label{sec21}

Let us assume that the function $f$ in \eqref{ODE} is sufficiently smooth, hence the initial-value problem has a unique smooth solution $y$, and we seek its approximation on an interval $[t_0,t_{\text{final}}]$. To this end, we apply a $k$-step LMM to \eqref{ODE} on a uniform grid $\{t_n\}$ to generate the sequence $y_n(h)$ according to \eqref{multistepdefn}. Here, $h:=t_{n+1}-t_n>0$ is the step size (or grid length), and $y_n(h)$ is supposed to approximate the exact solution at $t_n$, that is, $y_n(h)\approx y(t_n)$. We assume that the LMM is of order $p\ge 1$. 

The idea of classical RE is to take a suitable linear combination of two approximations, one generated on a coarser grid and one on a finer grid, to obtain a better approximation of the solution $y$ of \eqref{ODE}.  Here, we will only consider its simplest form and define
\begin{equation}\label{GREdef}
r_n(h):=\frac{2^p}{2^p-1}\cdot y_{2n}\left(\frac{h}{2}\right)-\frac{1}{2^p-1} \cdot y_n(h),
\end{equation}
that is, the coarser and finer grids have grid lengths $h$ and $h/2$, respectively.
Since the sequence $y_n(h)$ on the coarser grid  and the sequence $y_{2n}\left(\frac{h}{2}\right)$ on the finer grid are computed independently (their linear combination is formed only in the last step), we refer to this procedure as \textit{global} (or passive) \textit{extrapolation}, or, in short, LMM-GRE. 

\subsection{Convergence}\label{sec22}

\begin{lem}\label{convergencelemma} Under the above assumptions on the function $f$ in \eqref{ODE} and on the LMM, further, if  
the starting values $y_j(h)$ and $y_j\left(\frac{h}{2}\right)$ ($j=1,2,\ldots, k-1$) of the LMM are ${\cal{O}}(h^{p+1})$-close to the corresponding exact solution values, then the sequence $r_n(h)$ converges to the exact solution $y$ of \eqref{ODE}, and the order of convergence is at least $p+1$.
\end{lem}
\begin{proof}
The proof relies on the fact that---under the assumptions of the lemma---the global error $y_n(h)-y(t_n)$ of a LMM possesses an asymptotic expansion in $h$.
More precisely, according to, e.g., \cite[Section 6.3.4]{gautschi}, there exist a function $\mathbf{e}$ and a constant $C_{k,p}$ such that
\begin{equation}\label{globerr}
y_n(h)-y(t_n)=C_{k,p}\cdot h^p \cdot \mathbf{e}(t_n)+{\cal{O}}(h^{p+1})\quad\text{as }h\to 0^+,
\end{equation}
for any $n\in\mathbb{N}$ for which $t_n\in [t_0,t_{\text{final}}]$.
Here, 
the function $\mathbf{e}$ depends only on $f$ in \eqref{ODE} (and not on the chosen LMM), while the error constant $C_{k,p}$ depends only on the $k$-step LMM (and not on \eqref{ODE} or on $h$). Then, by applying \eqref{globerr} on a grid with grid length $h/2$ and focusing on the same (i.e., $h$-independent) grid point $t_n=t^*$, we have 
\begin{equation}\label{halfgloberr}
y_{2n}\left(\frac{h}{2}\right)-y(t^*)=C_{k,p}\cdot \left(\frac{h}{2}\right)^p \cdot \mathbf{e}(t^*)+{\cal{O}}(h^{p+1}).
\end{equation}
Combining \eqref{GREdef}--\eqref{halfgloberr}, we easily see that
$
r_n(h)-y(t^*)={\cal{O}}(h^{p+1})\quad\text{as }h\to 0^+.
$
\end{proof}

\subsection{Linear stability analysis}\label{sec23}

Let us now recall the definition of the region of absolute stability of a LMM---here, this
region will be denoted by $\mathcal{S}_\text{LMM}$. 
It is known (see \cite{hairerwanner} or \cite[Section 2.3]{optsubs}) that $\mathcal{S}_\text{LMM}$ can be characterized by the following boundedness condition. Let us fix some $h>0$ and $\lambda\in\mathbb{C}$ such that for $\mu:=h\lambda$ one has $\alpha_k-\mu\beta_k\ne 0$. Suppose that the LMM \eqref{multistepdefn}, with step size $h$ and starting values $y_0, y_1, \ldots, y_{k-1}$, applied to the usual scalar linear test equation 
\begin{equation}\label{Dahlquist}
y'(t)=\lambda y(t), \quad y(t_0)=y_0
\end{equation}
generates the sequence $y_n$ (${n\in\mathbb{N}}$). Then $\mu\in \mathcal{S}_\text{LMM}\subset\mathbb{C}$ if and only if the 
sequence $y_n$ is bounded for any choice of the starting values $ y_0, y_1, \ldots, y_{k-1}$.
\begin{rem}\label{orderreductionremark}
Considering the differential equation \eqref{Dahlquist}, if $\mu\in\mathbb{C}$ is chosen such that $\alpha_k-\mu\beta_k=0$, then the order of the recursion generated by the LMM becomes strictly less than $k$, hence the starting values $ y_0, y_1, \ldots, y_{k-1}$ could not be chosen arbitrarily (see also \cite[Remark 2.7]{optsubs}). 
\end{rem}

We define the region of absolute stability, ${\mathcal{S}}_\text{GRE}\subset\mathbb{C}$, of the combined LMM-GRE method \eqref{GREdef} analogously to that of the underlying LMM.
Let us apply \eqref{GREdef} to the scalar linear test equation \eqref{Dahlquist} with some $h>0$ and $\lambda\in\mathbb{C}$.
Then  ${\mathcal{S}}_\text{GRE}$ is defined to be\\

\textit{the set of numbers $\mu:=h\lambda$ for which the sequence $r_n(h)$ is bounded (in $n\in\mathbb{N}$) for any choice of the starting values of the sequence
$y_n(h)$ and for any choice of the starting values of the sequence $y_{m}\left(\frac{h}{2}\right)$, but excluding the  values of $\mu$ for which $\alpha_k-\mu\beta_k= 0$ or $\alpha_k-\frac{\mu}{2}\beta_k= 0$.}\\

 Now we can relate the stability region of the combined method to that of the underlying LMM as follows. For a set $S\subset\mathbb{C}$, we define $2\, S:=\{2z:z\in S\}$.
\begin{lem}\label{stabregionlemma} We have the inclusions $(\emph{i})\  {\mathcal{S}}_\text{\emph{LMM}}\cap \left(2\,{\mathcal{S}}_\text{\emph{LMM}}\right)\subseteq
{\mathcal{S}}_\text{\emph{GRE}}$,\ and $(\emph{ii})\   {\mathcal{S}}_\text{\emph{GRE}}\subseteq {\mathcal{S}}_\text{\emph{LMM}}$.
\end{lem}
\begin{proof}
Suppose that $h>0$ and $\lambda\in\mathbb{C}$ have been chosen such that $h\lambda\in{\mathcal{S}}_\text{{LMM}}\cap \left(2\,{\mathcal{S}}_\text{{LMM}}\right)$, and we apply the LMM-GRE method with this step size $h$ to \eqref{Dahlquist}  with this $\lambda$. Then both sequences $y_n(h)$ and $y_{m}\left(\frac{h}{2}\right)$ are bounded for any choice of their respective $k$ starting values. Hence the sequence $r_n(h)$, as their linear combination, is also bounded. This proves $(\mathrm{i})$.

To prove $(\mathrm{ii})$, let us choose $h>0$ and $\lambda\in\mathbb{C}$  such that $h\lambda\in{\mathcal{S}}_\text{\textrm{GRE}}$. Then the sequence $r_n(h)$ is bounded. By choosing every starting value $0$, we can have that the sequence $y_{m}\left(\frac{h}{2}\right)$ is identically $0$. Hence $r_n(h)=-\frac{1}{2^p-1} \cdot y_n(h)$, so the sequence $y_n(h)$ is also bounded. Therefore $h\lambda\in{\mathcal{S}}_\text{\textrm{LMM}}$.
\end{proof}
\begin{rem}
The reasoning in the above proof of $(\emph{ii})$ could not be applied to prove 
${\mathcal{S}}_\text{\emph{GRE}}\subseteq 2{\mathcal{S}}_\text{\emph{LMM}}$: although the boundedness of $r_n(h)$ implies (via a special choice of the starting values of the sequence $y_n(h)$) that the sequence $y_{2n}\left(\frac{h}{2}\right)$ is also bounded, this alone would be insufficient to guarantee the boundedness of the sequence 
$y_{m}\left(\frac{h}{2}\right)$ ($m\in\mathbb{N}$).
\end{rem}



To conclude this section, we give a sufficient condition for the stability regions ${\mathcal{S}}_\text{{LMM}}$ and ${\mathcal{S}}_\text{{GRE}}$ to coincide.
As it is well-know, all practically relevant LMMs are zero-stable \cite{suli}. 

\begin{lem}\label{convexlemma} Assume that the underlying LMM is zero-stable, and ${\mathcal{S}}_\text{\emph{LMM}}$ is convex. Then ${\mathcal{S}}_\text{\emph{GRE}}={\mathcal{S}}_\text{\emph{LMM}}$.
\end{lem}
\begin{proof} Zero-stability implies that $0\in{\mathcal{S}}_\text{{LMM}}$, so from the convexity of ${\mathcal{S}}_\text{{LMM}}$ we have that 
${\mathcal{S}}_\text{{LMM}}\subseteq 2{\mathcal{S}}_\text{{LMM}}$. But this means that ${\mathcal{S}}_\text{{LMM}}\cap (2{\mathcal{S}}_\text{{LMM}})={\mathcal{S}}_\text{{LMM}}$, so from Lemma \ref{stabregionlemma} we get that ${\mathcal{S}}_\text{{GRE}}={\mathcal{S}}_\text{{LMM}}$.
\end{proof}

\begin{rem} By analyzing the root-locus curve \cite{hairerwanner} of the underlying LMM as a parametric curve, it can be proved that ${\mathcal{S}}_\text{\emph{LMM}}$ is  convex, for example, for the Adams--Bashforth method with $k=2$ steps, or for the Adams--Moulton method with $k=2$ steps. However, for the Adams--Bashforth method with $k=3$ steps, ${\mathcal{S}}_\text{\emph{LMM}}$ is not convex.
\end{rem}

\subsection{$\text{BDF}k$-GRE methods}\label{sec24}

We obtain an efficient family of LMM-GRE methods if the underlying LMM is a $k$-step BDF method (referred to as a 
$\text{BDF}k$-method) with some $1\le k\le 6$  (recall that for zero-stability we need $k\le 6$). It is known that a $\text{BDF}k$-method has order $p=k$, see  \cite{hairerwanner}.

Suppose that the sequences $y_n(h)$ and $y_{2n}\left(\frac{h}{2}\right)$ in \eqref{GREdef} are generated by a $\text{BDF}k$-method, and the starting values for both sequences are $(k+1)^\text{st}$-order accurate. Then, due to Lemma \ref{convergencelemma}, \emph{the sequence $r_n(h)$ with $p:=k$ converges to the solution of \eqref{ODE} with order $k+1$}.

To measure the size of the region of absolute stability of the $\text{BDF}k$-GRE methods, one can invoke the concepts of $A$-stability and $A(\alpha)$-stability \cite{hairerwanner}.
It is easily seen that scaling the region of absolute stability of the underlying method ${\mathcal{S}}_\text{{LMM}}$ by a factor of $2$ preserves the $A(\alpha)$-stability angles (see \cite[Figure 1]{optsubs} for an illustration). Hence, due to Lemma \ref{stabregionlemma}, \emph{the  $\text{BDF}k$-GRE method has the same $A(\alpha)$-stability angle as that of the underlying $\text{BDF}k$-method}.


In Table \ref{tab:1}, we present the order of convergence and the $A(\alpha)$-stability angles for the $\text{BDF}k$-GRE methods. (For the \emph{exact} values of the angles $\alpha$, see, e.g., \cite[Table 1]{optsubs}.) The $\text{BDF}k$-GRE methods are particularly suitable for stiff problems.

\begin{table}[h]
\footnotesize
\caption{Convergence order and $A(\alpha)$-stability angles for the $\text{BDF}k$-GRE methods}
\label{tab:1}      
\centerline{\begin{tabular}{l|l|l}
\hline\noalign{\smallskip}
$k$ & order & $A(\alpha)$-stability angle\\
\noalign{\smallskip}\hline\noalign{\smallskip}
1 & $p=2$ & $90^\circ$, $A$-stable \\
2 & $p=3$ & $90^\circ$, $A$-stable \\
3 & $p=4$ & $86.032^\circ$ \\
4 & $p=5$ & $73.351^\circ$ \\
5 & $p=6$ & $51.839^\circ$ \\
6 & $p=7$ & $17.839^\circ$ \\
\noalign{\smallskip}\hline
\end{tabular}}
\end{table}






Notice, in particular, that the $\text{BDF}2$-GRE method is a $3^\text{rd}$-order $A$-stable method (recall that, due to the classical Dahlquist theorem, no $3^\text{rd}$-order $A$-stable LMM can exist). 

In terms of computational cost, due to the presence of the coarser and finer grids, the sequence $r_n(h)$ in \eqref{GREdef} corresponding to a LMM-GRE method is approximately three times as expensive to generate as the sequence $y_n(h)$ corresponding to the underlying LMM. However, the extra computing time is balanced by the higher order and $A(\alpha)$-stability; the $\text{BDF}5$-GRE method, for example, has order $6$, and its $A(\alpha)$-stability angle is approximately three times as large as the stability angle of the classical  $6^\text{th}$-order $\text{BDF}6$-method.


\section{Numerical experiments}\label{sec:numerics}

To verify the rate of convergence of LMM-GREs, we chose some benchmark problems,
including a Lotka--Volterra system
\[
y'_1(t)=0.1y_1(t)-0.3y_1(t)y_2(t), \quad y'_2(t)=0.5(y_1(t)-1)y_2(t)
\]
for $t\in[0,62]$ with initial condition $y(0)=(1,1)^\top$; or the mildly stiff {van der Pol} equation 
\[
y'_1(t)=y_2(t), \quad\quad y'_2(t)=2(1-y_1^2(t))y_2(t)-y_1(t)
\]
for $t\in[0,20]$ with initial condition $y(0)=(2,0)^\top$. 

\begin{figure}[ht!]
      \includegraphics[width=0.9\textwidth]{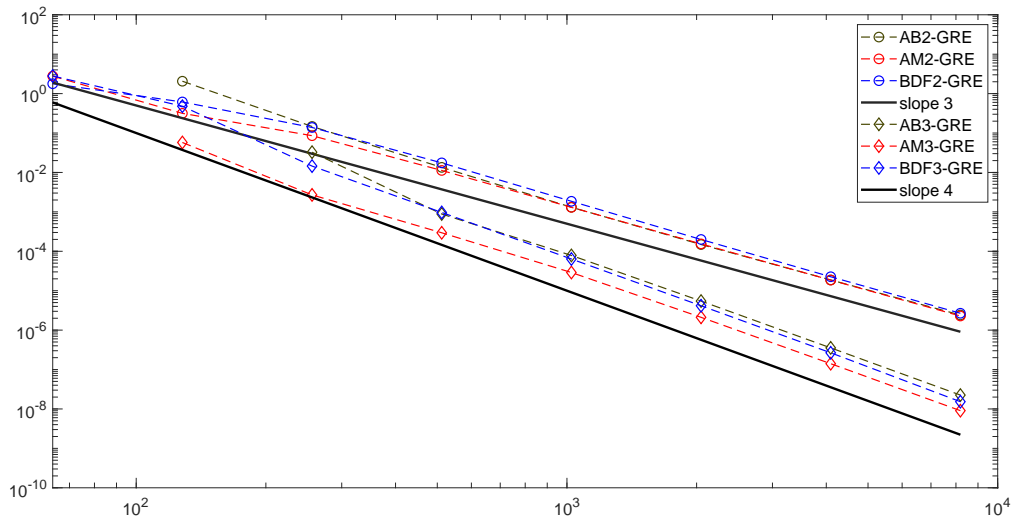}\caption{Results for the van der Pol equation for LMM-GREs, number of grid points versus the global error in maximum norm}\label{fig:1}
\end{figure}

As underlying LMMs, we considered the $2^\text{nd}$- and $3^\text{rd}$ order Adams--Bashforth (AB), Adams--Moulton (AM), and BDF methods. The AM methods were implemented in predictor-corrector style.  For starting methods, we chose the $2^\text{nd}$- and $3^\text{rd}$-order Ralston methods, having minimum error bounds \cite{ralston1962}. For the nonlinear algebraic equations arising in connection with implicit LMMs, we use MATLAB's \texttt{fsolve} command. Following \cite[Appendix A]{leveque}, the fine-grid solutions obtained by the classical $4^\text{th}$-order Runge--Kutta method with $2^{16}$ grid points are used to measure the global error in maximum norm and to estimate the corresponding convergence order. Table \ref{tab:2} and Figure \ref{fig:1} illustrate the expected order of convergence for all tested LMM-GREs. 

\begin{table}[ht!]
\begin{center}
\caption{The estimated order of convergence for the Lotka--Volterra system for different LMM-GREs with $64,128,\ldots,1024$ grid points}
\label{tab:2}
\begin{tabular}{ c|c|c|c|c|c}
\hline
$\text{AB}2$-GRE & $\text{AM}2$-GRE & $\text{BDF}2$-GRE & $\text{AB}3$-GRE &$\text{AM}3$-GRE & $\text{BDF}3$-GRE\\ 
 \hline
 3.7674 & 3.3545 & 3.2045 & 4.1864 & 3.7644 & 3.6254\\ 
 3.5761 & 3.2095 & 3.1703 & 3.9928 & 3.8903 & 3.7630\\ 
 3.1981 & 3.1068 & 3.1340 & 3.9873 & 3.9718 & 3.9411\\ 
 3.0297 & 3.0520 & 3.0807 & 3.9975 & 3.9856 & 3.9896\\
 \hline
\end{tabular}
\end{center}
\end{table}

\newpage

\subsection*{Acknowledgement}

The authors are indebted to an anonymous referee of the manuscript for their suggestions
that helped improving the presentation of the material, especially, for suggesting Lemma \ref{convexlemma} and its proof.\\

The project  ,,Application-domain specific highly reliable IT solutions'' has been implemented with the support provided from the National Research, Development and Innovation Fund of Hungary, financed under the Thematic Excellence Programme TKP2020-NKA-06 (National Challenges Subprogramme) funding scheme. I.~Fekete was supported by the J\'anos Bolyai Research Scholarship of the Hungarian Academy of Sciences, and also by the \'UNKP-21-5 New National Excellence Program of the Ministry for Innovation and Technology from the source of the National Research, Development and Innovation Fund.

\end{document}